\documentclass[11pt]{amsart2000} 
\usepackage{amssymb} 
\newcommand{\Q}{{\mathbb Q}} 
 
\newcommand{\N}{{\mathbb N}}
\newcommand{\Z}{{\mathbb Z}} 

\newtheorem{thm}{Theorem}

\begin{document}
\setlength{\unitlength}{0.01in}
\linethickness{0.01in}
\begin{center}
\begin{picture}(474,66)(0,0)
\multiput(0,66)(1,0){40}{\line(0,-1){24}}
\multiput(43,65)(1,-1){24}{\line(0,-1){40}}
\multiput(1,39)(1,-1){40}{\line(1,0){24}}
\multiput(70,2)(1,1){24}{\line(0,1){40}}
\multiput(72,0)(1,1){24}{\line(1,0){40}}
\multiput(97,66)(1,0){40}{\line(0,-1){40}}
\put(143,66){\makebox(0,0)[tl]{\footnotesize Proceedings of the Ninth Prague Topological Symposium}}
\put(143,50){\makebox(0,0)[tl]{\footnotesize Contributed papers from the symposium held in}}
\put(143,34){\makebox(0,0)[tl]{\footnotesize Prague, Czech Republic, August 19--25, 2001}}
\end{picture}
\end{center}
\vspace{0.25in}
\setcounter{page}{119}
\title[Universal minimal $S_\infty$-system]{The Cantor set of linear 
orders on ${\mathbb N}$ is the universal minimal $S_\infty$-system}
\author{Eli Glasner}
\address{Department of Mathematics\\
Tel Aviv University\\
Ramat Aviv\\
Israel}
\email{glasner@math.tau.ac.il}
\thanks{The author was an invited speaker at the Ninth Prague Topological Symposium.}
\thanks{This is a summary article. The results in this article will be 
treated fully in an article, written jointly with B. Weiss, to be 
published in Geometric and Functional Analysis (GAFA)}
\thanks{Eli Glasner,
{\em The Cantor set of linear orders on ${\mathbb N}$ is the universal 
minimal $S_\infty$-system},
Proceedings of the Ninth Prague Topological Symposium, (Prague, 2001),
pp.~119--123, Topology Atlas, Toronto, 2002}
\begin{abstract}
Each topological group $G$ admits a unique universal minimal dynamical 
system $(M(G),G)$. When $G$ is a non-compact locally compact group the 
phase space $M(G)$ of this universal system is non-metrizable.
There are however topological groups for which $M(G)$ is the trivial one 
point system (extremely amenable groups), as well as topological groups 
$G$ for which $M(G)$ is a metrizable space and for which there is
an explicit description of the dynamical system $(M(G),G)$.
One such group is the topological group $S_\infty$ of all permutations of 
the integers ${\mathbb Z}$, with the topology of pointwise convergence. 
We show that $(M(S_\infty),S_\infty)$ is a symbolic dynamical system
(hence in particular $M(S_\infty)$ is a Cantor set), and give a full 
description of all its symbolic factors.
Among other facts we show that $(M(G),G)$ (and hence also every minimal 
$S_\infty$) has the structure of a two-to-one group extension of proximal 
system and that it is uniquely ergodic.
\end{abstract}
\maketitle

This is a summary of a talk given at the Prague Topological Symposium
of 2001 in which I described results obtained in a joint paper with B. 
Weiss. The paper is going to appear soon in GAFA \cite{GW}.

Given a topological group $G$ and a compact Hausdorff space $X$, a 
dynamical system $(X,G)$ is a jointly continuous action of $G$ on $X$.
If $(Y,G)$ is a second dynamical system then a continuous onto map 
$\pi:(X,G)\to (Y,G)$ which intertwines the $G$ actions is called {\em a 
homomorphism\/}.
The dynamical system $(X,G)$ is {\em point transitive\/} if there exists a 
point $x_0\in X$ whose {\em orbit\/} $Gx_0$ is dense in $X$.
$(X,G)$ is {\em minimal\/} if every orbit is dense.
It can be easily shown that there exists a unique (up to isomorphism of 
dynamical $G$-systems) universal point transitive $G$-system 
$(\mathbf{L},G)$.
One way of presenting this universal object is via the Gelfand space of 
the $C^*$-algebra $\mathcal{L}_l(G)$ of left uniformly 
$\mathbb{C}$-valued continuous functions on $G$.
From the existence of $(\mathbf{L},G)$ one easily deduces the existence 
of a universal minimal dynamical system; i.e.\ a system $(M(G),G)$ such 
that for every minimal system $(X,G)$ there exists a homomorphism 
$\pi:(M(G),G)\to (X,G)$.
Ellis' theory shows that up to isomorphism this universal minimal 
dynamical system is unique, see e.g.\ \cite{E}.

The existence of uncountably many characters of the discrete group $\Z$
already shows that the phase space $M(\Z)$ is non-metrizable. In fact one
can show that $M(G)$ is non-metrizable whenever $G$ is non-compact locally
compact group.

A topological group $G$ has the {\em fixed point on compacta property\/} 
(f.p.c.) (or is {\em extremely amenable\/}) if whenever it acts 
continuously on a compact space, it has a fixed point. 
Thus the group $G$ has the f.p.c.\ property iff its universal minimal 
dynamical system is the trivial one point system.

A triple $(X,d,\mu)$, where $(X,d)$ is a metric space and $\mu$ a 
probability measure on $X$, is called an $mm${\em -space\/}.
For $A\subseteq X$, $\mu(A)\geq 1/2$, and $\epsilon>0$ let $A_\epsilon$ 
be the set of all points whose distance from $A$ is at most $\epsilon$.

A family of $mm$ spaces $(X_n,d_n,\mu_n)$ is called a {\em L\'evy 
family\/} if for every $\epsilon$, $\alpha_n(\epsilon)\to 0$, where
$\alpha(\epsilon)=1-\inf\{\mu(A_\epsilon) : A\subseteq X,\ \mu(A)\geq 
1/2\}$.
When a Polish group $(G,d)$ contains an increasing sequence of compact 
subgroups $\{G_n:n\in \N\}$ whose union is dense in $G$ and such that with 
respect to the corresponding sequence of Haar measures $\mu_n$, the family 
$(G_n,d,\mu_n)$ forms a L\'evy family, then $G$ is called a {\em L\'evy 
group\/}.

In \cite{GM} Gromov and Milman prove that every L\'evy group $G$ has the
f.p.c.\ property. Many of the examples presently known of extremely
amenable groups are obtained via this theorem. There are however other
methods of obtaining such groups. Here is a partial list:

\begin{enumerate}
\item
The unitary group $U(\infty)=\cup_{n=1}^\infty U(n)$ with the uniform 
operator topology (Gromov-Milman, \cite{GM}).
\item
The monothetic Polish group $L_m(I,S^1)$, consisting of all (classes) of 
measurable maps from the unit interval $I$ into the circle group $S^1$ 
with the topology of convergence in measure induced by, say, Lebesgue 
measure on $I$ (Glasner, \cite{G}; Furstenberg-Weiss).
More generally, $L_m(I,G)$, where $G$ is any locally compact amenable 
group (Pestov, \cite{P4}).
\item
The group of measurable automorphisms $\operatorname{Aut}(X,\mu)$ of a 
standard sigma-finite measure space $(X,\mu)$, with respect to the weak 
topology (Giordano-Pestov \cite{GP}).
\item
Using Ramsey's theorem, Pestov has shown that the group 
$\operatorname{Aut}({\Q,<})$, of order automorphism of the rational 
numbers with pointwise convergence topology, is extremely amenable, \cite{P2}.
\end{enumerate}

Thus, as we have seen, the universal minimal system $(M(G),G)$
corresponding to a non-compact $G$ is usually non-metrizable but can be,
in some cases, trivial. Are there non-compact topological groups for which
$M(G)$ is metrizable but non-trivial? The first such example was pointed
out by Pestov \cite{P2} who used claim 4 above to show that the universal
minimal dynamical system of the group $G$ of orientation-preserving
homeomorphisms of the circle coincides with the natural action of $G$ on
$S^1$.

In \cite{U} V. Uspenskij shows that the action of a topological group $G$
on its universal minimal system $M(G)$ is never $3$-transitive. As a
direct corollary he shows that for manifolds $X$ of dimension $>1$ (as
well as for $X=Q$, the Hilbert cube) the corresponding group $G$ of
orientation preserving homeomorphisms, $(M(G),G)$ does not coincide with
the natural action of $G$ on $X$.

Let $S_\infty$ be the group of all permutations of the integers $\Z$. With
respect to the topology of pointwise convergence on $\Z$, $S_\infty$ is a
Polish topological group. The subgroup $S_0\subset S_\infty$ consisting of
the permutations which fix all but a finite set in $\Z$ is an amenable
dense subgroup (being the union of an increasing sequence of finite
groups) and therefore $S_\infty$ is amenable as well.

In \cite{GM} Gromov and Milman conjectured, in view of the concentration 
of measure on $S_n$ with respect to Hamming distance, that $S_\infty$ has 
the f.p.c.\ property. In \cite{P2} and \cite{P3} V. Pestov has shown that, 
on the contrary, $S_\infty$ acts effectively on $M(S_\infty)$ and that, in 
fact, there is no Hausdorff topology making $S_0$ a topological group with 
the f.p.c.\ property. He as well as A. Kechris (in private communication) 
asked for explicit examples of $S_\infty$-minimal systems.

The main result of our work \cite{GW} is the fact that the universal
minimal system $(M(S_\infty),S_\infty)$ is a metrizable system, in fact a
system whose phase space is the Cantor set. We also give in this work an
explicit description of $(M(S_\infty),S_\infty)$ as a ``symbolic"
dynamical system and exhibit explicit formulas for all of its symbolic
factors. Let me now describe these results in more details.

For every integer $k\ge 2$ let
\begin{equation*}
\Z^k_*=\{(i_1,i_2,\dots,i_k)\in \Z^k:
i_1,i_2,\dots,i_k\ \text{ are distinct elements of $\Z$}\},
\end{equation*}
and set $\Omega^k=\{1,-1\}^{\Z^k_*}$.
Consider the dynamical system
$(\Omega^k,S_\infty)$, where for $\alpha\in S_\infty$ and $\omega\in 
\Omega^k$ 
we 
let
$$
(\alpha\omega)(i_1,i_2,\dots,i_k)=
\omega(\alpha^{-1}i_1,\alpha^{-1}i_2,\dots,\alpha^{-1}i_k).
$$

Let $\Omega^k_{alt}\subset\Omega^k$ consist of all the {\em alternating\/} 
configurations, that is those elements $\omega\in\Omega^k$ satisfying
$$
\omega(\sigma(i_1),\sigma(i_2),\dots,\sigma(i_k)) =
\operatorname{sgn}(\sigma)\omega(i_1,i_2,\dots,i_k),
$$
for all $\sigma\in S_k$ and $(i_1,i_2,\dots,i_k)\in \Z^k_*$.
Clearly $\Omega^k_{alt}$is a closed and $S_\infty$-invariant subset of 
$\Omega^k$.

A configuration $\omega\in\Omega^2$ {\em determines a linear order\/} on 
$\Z$ 
if 
it is alternating, and satisfies the conditions:
$$
\omega(m,n)=1\ \wedge\ \omega(n,l)=1 \quad \Rightarrow \quad 
\omega(m,l)=1.
$$
Let $<_\omega$ be the corresponding linear order on $\Z$, where 
$m<_\omega n$ iff $\omega(m,n)=1$.
Let $X=\Omega^2_{lo}$ be the subset of $\Omega^2$ consisting of all the 
configurations which determine a linear order.
The correspondence $\omega \longleftrightarrow\ <_\omega$ is a surjective 
bijection between $\Omega^2_{lo}$ and the collection of linear orders on 
$\Z$. Clearly $X$ is a closed $S_\infty$-invariant set and using Ramsey's 
theorem we shall show that $(X,S_\infty)$ is a minimal system.

Say that a configuration $\omega\in \Omega^3$ is {\em determined by a 
circular 
order\/} if there exists a sequence $\{z_m: m\in \Z\} \subset S^1$ with 
$m\ne n \ \Rightarrow\ z_m\ne z_n$ such that: $\omega(l,m,n)=1$ for 
$(l,m,n)\in 
\Z^3_*$ iff the directed arc in $S^1$ defined by the ordered triple 
$(z_l,z_m,z_n)$ is oriented in the positive direction. 
Let $Y=\Omega^3_c\subset \Omega^3_{alt}$ denote the collection of all the 
configurations in $\Omega^3$ which are determined by a circular order.
It follows that the set $Y=\Omega^3_c$ is closed and invariant and using 
Ramsey's theorem one can show that it is minimal.

If we go now to $\Omega^4_{alt}$, can one find a sequence of points 
$\{z_n\}$ 
on the sphere $S^2$ in general position such that the tetrahedron defined 
by any four points $z_{n_1},z_{n_2},z_{n_3},z_{n_4}$ has positive 
orientation when $n_1< n_2 < n_3 < n_4$?

Starting with any sequence $\{z_n\}\subset S^2$ in general position one 
can use Ramsey's theorem to find a subsequence with the required property. 
Another way to see this is to use the `moment curve'
$$
t\mapsto (t,t^2,t^3).
$$
Again it turns out that the orbit closure in $\Omega^4_{alt}$ which is 
determined by such a sequence forms a minimal dynamical system.

It now seems as if going up to $\Omega^k_{alt}$ with larger and larger 
$k$'s 
we encounter more and more complicated minimal systems.
However, as we show in \cite{GW}, this is not the case and the entire 
story is already encoded in the simplest symbolic dynamical system 
$\Omega^2_{lo}$.

\begin{thm}
$\Omega^2_{lo}$ is the universal minimal $S_\infty$-system.
\end{thm}

The fact that the topology on $S_\infty$ is zero-dimensional,
and in fact given by a sequence of clopen subgroups,
enables us to reduce this theorem to the following one.

\begin{thm}
Every minimal subsystem $\Sigma$ of the system $(\Omega^k,S_\infty)$ is
a factor of the minimal system $(\Omega^2_{lo},S_\infty)$.
\end{thm}

Finally let me mention two more facts concerning the system 
$(M(S_\infty),S_\infty)$.

\begin{thm}
The universal minimal system $(\Omega^2_{lo},S_\infty)$ has the structure
of a two-to-one group extension of a proximal system.
\end{thm}

\begin{thm}
The universal minimal system $(\Omega^2_{lo},S_\infty)$ is uniquely 
ergodic
and therefore so is every minimal $S_\infty$-system.
\end{thm}

\providecommand{\bysame}{\leavevmode\hbox to3em{\hrulefill}\thinspace}
\providecommand{\MR}{\relax\ifhmode\unskip\space\fi MR }
\providecommand{\MRhref}[2]{%
  \href{http://www.ams.org/mathscinet-getitem?mr=#1}{#2}
}
\providecommand{\href}[2]{#2}

\end{document}